\documentclass{scrartcl}
\usepackage{amsmath,amssymb,amsthm}
\usepackage{xcolor}
\usepackage{dsfont}
\usepackage{mathtools}
\usepackage[numbers]{natbib}
\usepackage{enumerate}

\usepackage{thmtools}
\usepackage{hyperref, cleveref}
\usepackage{mathdots}

\theoremstyle{plain}
\newtheorem{thm}{Theorem}[section]
\newtheorem{cor}[thm]{Corollary} 
\newtheorem{lem}[thm]{Lemma}
\newtheorem{conj}[thm]{Conjecture}
\newtheorem{prop}[thm]{Proposition}
\theoremstyle{definition}
\newtheorem{defn}[thm]{Definition}
\newtheorem{obs}[thm]{Observation}

\DeclareMathOperator{\inv}{inv}
\DeclareMathOperator{\rk}{rank}
\DeclareMathOperator{\tmr}{tmr}

\newcommand\F{\mathbb{F}}
\newcommand\dij{D_1\rightarrow D_2}
\DeclareMathAlphabet{\mathmybb}{U}{bbold}{m}{n}
\newcommand{\1}{\mathmybb{1}}

\title{A case of the dijoin conjecture on inverting oriented graphs}
\author{Natalie Behague\thanks{School of Mathematical Sciences, Dublin City University, Dublin, Ireland. \\Email: \href{mailto:natalie.behague@dcu.ie}{natalie.behague@dcu.ie}} \and Patrick Gaudart-Wifling\thanks{School of Mathematical Sciences, Queen Mary University of London, London, UK.\\Email: \href{mailto:p.gaudart-wifling@qmul.ac.uk}{p.gaudart-wifling@qmul.ac.uk}}}

\date{September 2025}

\begin{document}

\maketitle
\begin{abstract}
    For an oriented graph $D$, the inversion of $X\subseteq V(D)$ in $D$ is the graph obtained by reversing the orientation of all arcs with both ends in $X$. The inversion number $\inv(D)$ is the minimum number of inversions needed to obtain an acyclic oriented graph.

    We show that the dijoin conjecture of Bang-Jensen, da Silva and Havet, that $\inv(\dij)=\inv(D_1)+\inv(D_2)$, is true in the case where $\inv(D_1)=2$ and $\inv(D_2)$ is even. We also characterise the cases $\inv(D_1)=2$ and $\inv(D_2)$ odd, for which the conjecture does and does not hold. We then go on to show a similar result for t-joins, in doing so we prove a conjecture of Alon, Powierski, Savery, Scott and Wilmer. Our proofs build on the idea of tournament minimum rank, introduced by Behague, Johnston, Morrison and Ogden.
\end{abstract}

\section{Introduction}
For an oriented graph $D$ and a set $X\subseteq V(D)$, the \emph{inversion} of $X$ in $D$ is the graph obtained by reversing the orientation of every arc with both ends in $X$. The inversion of a collection $X_1, X_2, \dots, X_m$, is the graph obtained by successively inverting each set. Note that the order in which this is done does not matter. $X_1, X_2, \dots,X_m$ is called a \emph{decycling family} for $D$ if its inversion contains no directed cycle. 

Belkhechine \citep{belkhechine2009} introduced the inversion number $\inv(D)$, defined as the smallest integer $m$ for which an $m$-decycling family, $X_1,X_2,\dots,X_m$, of $D$ exists. Initial results on inversions were then found by Belkhechine, Bouaziz, Boudabbous and Pouzet \citep{belkhechine2010}, although they were mainly concerned with the quantity $\inv(n)$, which is the maximum inversion number for an oriented graph of order $n$. 

More recently, Bang-Jensen, da Silva and Havet \citep{bang2022} examined the inversion number for dijoins. For oriented graphs $D_1$ and $D_2$, their \emph{dijoin}, denoted $\dij$ is defined by taking the disjoint union of $D_1$ and $D_2$ and adding an arc from each $u\in V(D_1)$ to each $v\in V(D_2)$. They conjectured that the obvious upper bound for the inversion number of a dijoin is actually an equality: ${\inv(\dij)=\inv(D_1)+\inv(D_2)}$. They proved their conjecture in the case $\inv(D_1)+\inv(D_2)\leq3$, and the authors of \citep{alon2024} proved it for $\inv(D_1)=\inv(D_2)=2$.

However, the authors of \citep{alon2024} and \citep{aubian2025} have since disproven the conjecture in general. In particular, \citep{aubian2025} showed that for each odd $k\geq3$, there exists a tournament $D_1$ with $\inv(D_1)=k$, such that for any other oriented graph $D_2$ with $\inv(D_2)\geq1$, we have $\inv(\dij)\leq\inv(D_1)+\inv(D_2)-1$. 

Two simultaneous papers \citep{behague2025,wang2024} have examined the quantity $\inv(\dij)$ further. Behague, Johnston, Morrison and Ogden \citep{behague2025} made an initial step in finding a lower bound, showing that for oriented graphs $D_1$ and $D_2$ with $\inv(D_1)=\inv(D_2)\geq1$, we have $\inv(\dij)>\inv(D_1)$. Wang, Yang and Lu \citep{wang2024} showed that the dijoin conjecture does still hold in the case where $\inv(D_1)=1$ and $\inv(D_2)$ is even. 

We combine and build on the methods used in these two papers. Before we state our main result, let us introduce a key definition.
\begin{defn}\label{def:inv_rank_distinct}
An oriented graph $D$ is called \emph{overfull} if $\inv(\overrightarrow{C_3} \rightarrow D) = inv(D)$, where $\overrightarrow{C_3}$ denotes the directed cycle on 3 vertices.
\end{defn}

Our main theorem is the following:
\begin{restatable}{thm}{dijoin}\label{thm:dijoin}
    Let $D_1$ and $D_2$ be oriented graphs such that $\inv(D_1)\in\{1,2\}$, then either:
    \begin{enumerate}[i.]
        \item $\inv(\dij)=\inv(D_1) + \inv(D_2)$; or
        \item $\inv(\dij)=\inv(D_1) + \inv(D_2) -1$.
    \end{enumerate}
    Moreover, case (ii) holds if and only if $D_2$ is overfull, in which case $\inv(D_2)$ is odd.
\end{restatable}

\Cref{thm:dijoin} proves that the dijoin conjecture holds for $\inv(D_1)=2$ and $\inv(D_2)$ even. Note that the $\inv(D_1) = 1$ case of \Cref{thm:dijoin} follows immediately from the work of Wang, Yang and Lu~\citep[Corollary~1.10]{wang2024}.

 Let $[D_1,D_2,\dots,D_t]$ be the oriented graph $(\ldots ((D_1 \rightarrow D_2) \rightarrow D_3)\rightarrow\ldots \rightarrow D_t)$ obtained by taking successive dijoins, the \emph{t-join} (note that the order the dijoins are carried out does not matter). Alon, Powierski, Savery, Scott and Wilmer made the following conjecture \citep{alon2024}.

\begin{conj}[{\citep[Conjecture~9]{alon2024}}]\label{conj:alon}
    Let $t\in\mathbb{N}$ and let $D_1,\dots,D_t$ be oriented graphs satisfying $\inv(D_i)\leq2$ for all $i$. Then
    \[
    \inv([D_1,\dots,D_t])=\sum_{i=1}^t\inv(D_i).
    \]
\end{conj}

Wang, Yang and Lu proved the following result, which resolves some cases of the conjecture, as well as some cases not covered by the conjecture.
\begin{thm}[{\citep[Corollary~1.10]{wang2024}}] \label{thm:Wang k-join}
    Let $t\geq2$ and $D_1,D_2,\dots,D_t$ be oriented graphs. Suppose that there is some $j\in[t]$ such that $\inv(D_i)=1$ for all $i\in[t]\backslash\{j\}$ and $\inv(D_j)\geq1$. Then
    \[
    \inv([D_1,D_2,\dots,D_t])=\begin{cases}
        \sum_{i=1}^t\inv(D_i)-1 & \text{if $D_j$ is overfull},
        \\ \sum_{i=1}^t\inv(D_i) & \text{otherwise.}
    \end{cases}
    \]
\end{thm}

We are able to settle~\Cref{conj:alon}. In fact, we prove the following slightly stronger result:

\begin{restatable}{thm}{kjoin}\label{thm:k-join}
    Let $t\geq2$ and $D_1,D_2,\dots,D_t$ be oriented graphs. Suppose that there is some $j\in[t]$ such that $\inv(D_i)\leq2$ for all $i\in[t]\backslash\{j\}$ and $\inv(D_j)\geq0$. Then
    \[
    \inv([D_1,D_2,\dots,D_t])=\begin{cases}
        \sum_{i=1}^t\inv(D_i)-1 & \text{if $D_j$ is overfull and} \\  & \qquad \text{$\inv(D_i) \ne 0$ for some $i \ne j$},
        \\ \sum_{i=1}^t\inv(D_i) & \text{otherwise.}
    \end{cases}
    \]
\end{restatable}

\subsection{Paper structure}
In \Cref{sec:prelims} we will introduce key definitions, most notably that of tournament minimum rank, and some related results. In \Cref{sec:tmr} we prove key lemmas on the tournament minimum rank. In \Cref{sec:inv tourney} we apply these lemmas to prove versions of \Cref{thm:dijoin} and \Cref{thm:k-join} for tournaments, and then in \Cref{sec:inv all} we extend these results to all oriented graphs. We conclude with some directions and ideas for future research on this topic in \Cref{sec:Conjectures}.

\subsection{Notation}
    Given a transitive tournament $T$, we use the sign `$<$' to express the order of its vertices, that is $x<y$ if $xy\in A(T)$. Then $x<y<z$ means $xy,xz,yz\in A(T)$. For sets $A,B\subseteq V(T)$, $A<B$ if $x<y$ for all $x\in A$ and $y\in B$.

\section{Preliminaries}\label{sec:prelims}

In this paper the oriented graphs will always be simple, that is, with no loops and at most one arc between any two vertices. A \emph{tournament} is an oriented graph with exactly one arc between any two vertices. 

\begin{defn}
    Let $D$ be a tournament with $V(D)=\{v_1,v_2,\dots,v_n\}$. We say that an $n\times n$ symmetric matrix $M$ with entries in $\F_2$ is a \emph{decycling matrix} for $D$, if reversing the orientation of all arcs $v_iv_j$ where $m_{ij}=1$, results in a tournament with no directed cycles.     
\end{defn}
\begin{defn}
Let $D$ be a tournament. The \emph{tournament minimum rank} of $D$, denoted $\tmr(D)$, is the smallest rank of a decycling matrix for $D$.
\end{defn}

\begin{defn}
    For a tournament $D$, a decycling family $X_1, X_2,\dots,X_m$, and a vertex $v\in V(D)$, we define the \emph{characteristic vector} $\mathbf{v}\in\F_2^m$ of $v$, where the $i^{\emph{th}}$ element of $\mathbf{v}$ is 1 if and only if $v\in X_i$.
\end{defn}

Notice that for two vertices $v,v'\in V(D)$, $\mathbf{v\cdot v'}=1$ if and only if their arc has opposite orientations, before and after the decycling family. The \emph{gram matrix} of a set of vectors $v_1, \ldots, v_n$ has $i,j$th entry $v_i \cdot v_j$.  Therefore, the gram matrix of a decycling family's characteristic vectors is a decycling matrix for $D$. 

We can also pass from a decycling matrix to a set of characteristic vectors, due to a result of Buchanan, Purcell and Rombach \citep{buchanan2022}:

\begin{lem}[{\citep[Theorem~4.6]{buchanan2022}}]\label{lem: buchanan}
    Let $A$ be an $n\times n$ symmetric matrix over $\F_2$ of rank $k$. Then either,
    \begin{enumerate}[i.]
        \item $A=XX^T$, where $X$ is an $n\times k$ matrix over $\F_2$ of rank $k$; or
        \item $A=XX^T$, where $X$ is an $n\times (k+1)$ matrix over $\F_2$ of rank $k$ and $k$ is even. 
    \end{enumerate}
\end{lem}

By taking $A$ in this lemma to be a decycling matrix for some tournament $D$, we see that the rows of $X$ will be the characteristic vectors of a decycling family for $D$ ($A$ is their gram matrix). 

The relationship between tournament minimum rank and the inversion number was given by Behague, Johnston, Morrison and Ogden \citep{behague2025} in the following result: 
\begin{cor}[{\citep[Corollary~3.2]{behague2025}}]\label{cor:BJMO main}
    Let $D$ be a tournament. Then either
    \begin{enumerate}
        \item $\inv(D)=\tmr(D)$, or
        \item $\inv(D)=\tmr(D)+1$, in which case $\tmr(D)$ is even.
    \end{enumerate}
    Moreover, if $D$ is not transitive, then $\inv(D)=\tmr(D)+1$ if and only if every decycling matrix for $D$ with minimum rank has every diagonal entry equal to zero.  
\end{cor}

The following lemma illustrates why tournament minimum rank is particularly useful when discussing the dijoin conjecture. The proof is inspired by the construction of counterexamples in \citep{aubian2025}. 

\begin{lem}\label{lem: inv + tmr}
    Let $D_1$ and $D_2$ be tournaments, with $\inv(D_1)\geq1$. Then $\inv(\dij)\leq\inv(D_1)+\tmr(D_2).$
\end{lem}

\begin{proof}
    Let $\ell =\inv(D_1)$ and $k=\tmr(D_2)$. We will construct an $(\ell+k)$-decycling family for $\dij$.

    Let $X_1,X_2,\dots,X_\ell$ be a decycling family for $D_1$, with characteristic vectors $\mathbf{v}_1,\mathbf{v}_2,\dots,\mathbf{v}_m$.

    If $\inv(D_2)=\tmr(D_2)=k$, we can take a $k$-decycling family $X'_1,X'_2,\dots,X'_k$ for $D_2$ and we get the decycling family for $\dij$: 
    \[
    X_1,X_2,\dots,X_\ell,X'_1,X'_2,\dots,X'_k.
    \]

    Suppose instead that $\inv(D_2)=\tmr(D_2)+1=k+1$. Then, by \Cref{cor:BJMO main},  $k$ is even and every minimum rank decycling matrix for $D_2$ has diagonal entries equal to zero. Take such a matrix $B$ then, by \Cref{lem: buchanan}, we can find a rank $k$, $n\times(k+1)$ matrix $Y$ such that $YY^T=B$. We therefore have that the rows of $Y$ are characteristic vectors for a $(k+1)$-decycling family of $D_2$, $X'_1,X'_2,\dots,X'_{k+1}$.
    
    Call these vectors $\mathbf{y}_1,\mathbf{y}_2,\dots,\mathbf{y}_n$. Notice that the diagonal entries of $B$ are $\mathbf{y}_i\cdot\mathbf{y}_i$ and so, since they are all zero, we have for all $i$:
    \[
    \mathbf{1}\cdot\mathbf{y}_i=\mathbf{y}_i\cdot\mathbf{y}_i=0
    \]
    where $\mathbf{1}=(1,1,\dots,1)$. Note also that, clearly, $\mathbf{y}_i\cdot\mathbf{0}=0$. 

    We therefore claim that we have the following $(\ell+k)$-decycling family for $\dij$:
    \[
    X_1'\cup X_1,\; X_2'\cup X_1,\; \dots,\; X_k'\cup X_1,\; X'_{k+1}\cup X_1,\; X_2,\; X_3, \dots,\; X_\ell.
    \]

    Notice that the characteristic vectors for vertices in $D_1$ are $(\mathbf{1},\mathbf{v}_i)$ or $(\mathbf{0},\mathbf{v}_i)$. The characteristic vectors for vertices in $D_2$ are $(\mathbf{y}_i,\mathbf{0})$. 

    Therefore, the gram matrix for the whole set of characteristic vectors is $\begin{psmallmatrix}
        A & 0\\ 0 & B
    \end{psmallmatrix}$, where $A$ and $B$ are decycling matrices for $D_1$ and $D_2$ respectively. This is a decycling matrix for $\dij$, and so
    $\inv(\dij) \leq \ell+k=\inv(D_1)+\tmr(D_2)$.

\end{proof}

We can restate the property of being `overfull' (\Cref{def:inv_rank_distinct}) in terms of tournament minimum rank. To do so, we use two results of Wang, Yang and Lu, together with a result of Behague, Johnston, Morrison and Ogden.

\begin{lem}[{\citep[Lemma~2.6]{wang2024}}]\label{lem: wang1}
    Let $D$ be an oriented graph. If $\inv(D)=k\ge 3$, then $D$ is overfull if and only if  $k$ is odd and there exists a $k$-decycling family for $D$ whose decycling matrix has rank $k-1$ and all diagonal entries equal to zero. 
\end{lem}
\begin{lem}[{$k=1$: \citep[Theorem~2]{behague2025} and $k=0,2$: \citep[Theorem~1.7]{wang2024}}]\label{lem: small_inv_not_overfull}
     Let $D$ be an oriented graph. If $\inv(D)\le 2$ then $D$ is not overfull. 
\end{lem}
Combining this statement with \Cref{cor:BJMO main} we get the following corollary:

\begin{cor}\label{cor: YWL to tmr}
    Let $D$ be a tournament. Then $D$ is overfull if and only if $\inv(D)=\tmr(D)+1$, in which case $\inv(D)$ is odd. Equivalently, $D$ is not overfull if and only if  $\inv(D)=\tmr(D)$.
\end{cor}
\begin{proof}
    First suppose $\inv(D)=\tmr(D)+1$. Note that this trivially implies $inv(D) \ne 1$. By \Cref{cor:BJMO main}, $\inv(D)$ is odd and its minimum rank decycling matrices have all diagonal entries equal to zero. Take such a matrix, then by \Cref{lem: buchanan} we can find a decycling family of length $\inv(D)$ that generates it.  Since $\inv(D)$ is odd and $\ne 1$, we can apply  \Cref{lem: wang1} to conclude $\inv(\overrightarrow{C_3}\rightarrow D)=\inv(D)$.
    
    For the converse, suppose that $\inv(\overrightarrow{C_3}\rightarrow D)=\inv(D)$. By \Cref{lem: small_inv_not_overfull}, we have $\inv(D) \ge 3$. Thus by \Cref{lem: wang1}, we see that $\tmr(D)\le\inv(D)-1$. By~\Cref{cor:BJMO main},  $\tmr(D)\ge\inv(D)-1$ and so we have equality, as required. 
\end{proof}

\section{Results on Tournament Minimum Rank}\label{sec:tmr}

The aim of this section is to prove analogues of Theorems \ref{thm:dijoin} and \ref{thm:k-join}, for the tournament minimum rank rather than the inversion number. A \emph{principal submatrix} of a square matrix is a submatrix obtained by deleting rows and columns indexed by the same set. We need the following lemma: 

\begin{lem}\label{lem: submatrix}
    Let M be a symmetric $n\times n$ matrix with entries in $\F_2$. If $\rk(M)\geq2$ then $M$ contains a principal $2\times2$ submatrix with rank 2.
\end{lem}
\begin{proof}
    Suppose for a contradiction that $M$ has no such principal submatrix. Order its columns, and respective rows, such that we have:
    \[
    \begin{pmatrix}
        A & C \\
        C^T & B
    \end{pmatrix}
    \]
    where the diagonal entries in $A$ are all 1s and the diagonal entries in $B$ are all 0s.

Note that $B$ is either empty or must consist entirely of 0s, else we would find a principal $\begin{psmallmatrix}
        0 & 1 \\ 1 & 0
        \end{psmallmatrix}$.
        
        If $A$ is empty, then $M = B$ and $\rk(M) = 0$. 
    Thus $A$ is non-empty.  $A$ must consist entirely of 1s, else we would find a principal $\begin{psmallmatrix}
        1 & 0 \\ 0 & 1
    \end{psmallmatrix}$. Similarly $C$ must consist entirely of 0s, else we would find a principal $\begin{psmallmatrix}
            1 & 1 \\ 1 & 0
        \end{psmallmatrix}$. Therefore, $\rk(M)<2$.

\end{proof}

We can now give our first result for tournament minimum rank. 

\begin{lem}\label{lem: tmr + tmr}
    Let $D_1,D_2$ be tournaments such that $\tmr(D_1)=2$ and $\tmr(D_2)=k$. Then $\tmr(\dij)=k+2.$
\end{lem}
\begin{proof}
    We know $\tmr(\dij)\leq k+2$, since we can take minimum rank decycling matrices $A$ and $B$, for $D_1$ and $D_2$ respectively, and $\begin{psmallmatrix}
        A & 0 \\ 0 & B
    \end{psmallmatrix}$ is a decycling matrix for $\dij$. 

    Suppose for a contradiction that $\tmr(\dij)\leq k+1$. We will find a decycling matrix for $D_2$ with rank $k-1$. 

    Let $M$ be a decycling matrix for $D_1 \rightarrow D_2$ with rank $\leq k+1$, such that the ordering for the vertices in the resulting transitive tournament is:
    \[
    V_0<u_1<V_1<u_2<V_2<\dots<u_n<V_n
    \]
    where $\{u_1,u_2,\dots,u_n\}=V(D_1)$ and $V_0,V_1,\dots,V_n$ is a partition of $V(D_2)$ (with possibly empty sets). 

    We then order the rows/columns of $M$ according to the vertex ordering $u_1,u_2,\dots,u_n,V_0,V_1,\dots,V_n$, so that $M$ has the form:
    \[
    \begin{pmatrix}
        A & C \\
         C^T & B 
    \end{pmatrix}
    \]
    
    where, $A$ and $B$ are decycling matrices for $D_1$ and $D_2$ respectively. Notice that the entries of $C$ increase down each column and decrease along each row (so that the 1s form the shape of a staircase in the bottom left).

    By \Cref{lem: submatrix}, since $\rk(A)\geq2$, we can find a principal $2\times 2$ submatrix $A'$ of $A$ with rank 2, corresponding to vertices $u_i$ and $u_j$. We then combine the sets of $V(D_2)$ further based on whether or not their arcs with $u_i$ and $u_j$ are flipped by the inversions: 
    \begin{itemize}
        \item $P=V_0\cup\dots\cup V_{i-1}$
        \item $Q=V_i\cup\dots\cup V_{j-1}$
        \item $R=V_j\cup\dots\cup V_n$
    \end{itemize}
    Now consider the submatrix of $M$ given by $u_i,u_j,P,Q,R$:

\[    
\begin{pmatrix}
    
    a_{ii} & a_{ij} & 1 \dots 1 & 0 \dots 0 & 0 \dots 0 \\
    a_{ji} & a_{jj} & 1 \dots 1 & 1 \dots 1 & 0 \dots 0 \\
    1 & 1 & b_{11}\dots\null & \dots & \text{ }\dots b_{1n}\\
    \vdots & \vdots & \vdots \text{ } \ddots & & \iddots \text{}\vdots\\
    1 & 1 \\
    0 & 1 \\
    \vdots & \vdots & \vdots &  & \vdots \\
    0 & 1 \\
    0 & 0 \\
    \vdots & \vdots & \vdots \text{ } \iddots & & \ddots \text{ } \vdots \\
    0 & 0 & b_{1n}\dots\null & \dots & \text{ }\dots b_{nn}
    
\end{pmatrix}.
\]

Since $\rk(A')=2$, we can add some combination of columns 1 and 2 to the third column to make entries $M_{1,3}, M_{2,3}$ both $0$. Since $A'$ is symmetric (as $A$ is), we can add the respective combination of rows 1 and 2 to the third row to make entries $M_{3,1}, M_{3,2}$ both $0$. Continuing in this way for subsequent rows and columns -- adding some combination of columns 1 and 2 to the $i$th column and the respective combination of rows 1 and 2 to the the $i$th row -- we end up with a matrix of the form
\[
\begin{pmatrix}
    A' & 0 \\ 0 & B'
\end{pmatrix}.
\]
where $B'$ is symmetric.
This matrix still has rank at most $k+1$ and so $\rk(B')\leq k-1$. It remains to check that $B'$ is still a decycling matrix for $D_2$. 

The order of the vertices of $D_2$ after the inversions of $B$ is $P<Q<R$. So we can consider $B$ as a block matrix: 
\[
\begin{pmatrix}
B_{PP} & B_{PQ} & B_{PR} \\
B_{QP} & B_{QQ} & B_{QR} \\
B_{RP} & B_{RQ} & B_{RR} \\
\end{pmatrix}
\]

The row/column additions do not change any of the blocks corresponding to arcs with a vertex in $R$, so these arcs have the same orientation after $B$ or $B'$.

The entries in the block $B_{PP}$ each have the same combination added to them and so these entries are either all the same, or all changed. Therefore, the arcs with both vertices in $P$ are either all opposite or all the same after $B$ and $B'$, either way $P$ remains acyclic.

By the same argument, $Q$ remains acyclic. 

Finally, all entries in the blocks $B_{QP}$ and $B_{PQ}$  receive the same combination of row/column additions, so these entries are either all the same, or all changed. So the arcs between $P$ and $Q$ are either all opposite or all the same after $B$ and $B'$. 

So $B'$ is a decycling matrix. Its final ordering is either $P<Q<R$ or $Q<P<R$, where $P, Q$ and $R$ are each acyclic.

\end{proof}

Notice that, if $\tmr(D_1)\geq2$, \Cref{lem: submatrix} allows us to follow the same argument and we get a lower bound:

\begin{cor}\label{cor:tmr bound}
    Let $D_1,D_2$ be tournaments such that $\tmr(D_1)\geq2$ and $\tmr(D_2)=k$. Then $\tmr(\dij)\geq k+2.$
\end{cor}

We would now like to get a result for the tournament minimum rank of $n$-joins. 
The first step is to show that the additivity property in \Cref{lem: tmr + tmr} also holds when $\tmr(D_1)=1$. 
In fact, this follows from the results of Wang, Yang and Lu \citep{wang2024}, although they did not use tournament minimum rank in their arguments.

\begin{lem}[{\citep[Lemma 2.7]{wang2024}}]\label{lem: wangb}
    Let $D_1$ and $D_2$ be oriented graphs such that $\inv(D_1)=1$. Then $\inv(\dij)=\inv(\overrightarrow{C_3}\rightarrow D_2)$, and $\inv(D_2\rightarrow D_1)=\inv(D_2\rightarrow \overrightarrow{C_3})$. 
\end{lem}

Wang, Yang and Lu \citep{wang2024} also showed the following:

\begin{lem}[{\citep[Corollary~1.10]{wang2024}}]\label{lem: wang2}
    Let $D_1,D_2$ be oriented graphs such that $\inv(D_1)=1$. Then $\inv(\overrightarrow{C_3}\rightarrow(\dij))=1+\inv(\dij)$. In particular, $D_1 \rightarrow D_2$ is not overfull.
\end{lem}

We therefore get the desired property, similar to \Cref{lem: tmr + tmr}:

\begin{cor}\label{cor:tmr(1->k)}
    Let $D_1,D_2$ be tournaments such that $\tmr(D_1)=1$ and $\tmr(D_2)=k$. Then $\tmr(\dij)=k+1.$
\end{cor}
\begin{proof}
    From \Cref{lem: wang2}, $\dij$ is not overfull. Thus by \Cref{cor: YWL to tmr}, we have $\inv(\dij)= \tmr(\dij)$.    
    If $\inv(\dij)=k$, then by \Cref{lem: wangb}, $\inv(\overrightarrow{C_3}\rightarrow D_2)=k$ and so $D_2$ is overfull. By  \Cref{cor: YWL to tmr}, we get $\tmr(D_2)=k-1$, a contradiction.  Therefore $\tmr(\dij)=\inv(\dij)=k+1$.
\end{proof} 

Using the additive property shown in \Cref{lem: tmr + tmr} and \Cref{cor:tmr(1->k)}, we get the tournament minimum rank version of \Cref{thm:k-join}. 

\begin{lem}\label{lem:tmr sum}
    Let $t\geq2$ and $D_1,D_2,\dots,D_t$ be tournaments. Suppose that there is $j\in[t]$ such that $\tmr(D_j)\geq1$ and $\tmr(D_i)=1\text{ or }2$ for all $i\in[t]\backslash\{j\}$. Then,
    \[
    \tmr([D_1,D_2,\dots,D_t])=\sum_{i=1}^t\tmr(D_i).
    \]
\end{lem}

In order to prove this, we need the following proposition:

\begin{prop}
    Let $D_1$ and $D_2$ be tournaments such that $\tmr(D_1)=1\text{ or }2$ and $\tmr(D_2)=k$. Then $\tmr(D_2\rightarrow D_1)=\tmr(\dij)=\tmr(D_1)+\tmr(D_2)$.
\end{prop}
\begin{proof}
    Denote by $D^-$ the tournament obtained by reversing all arcs of $D$, and notice that $D$ and $D^-$ have the same set of decycling matrices (the resulting acyclic tournaments will have opposite orientations), so $\tmr(D)=\tmr(D^-)$. We therefore have
    \begin{align*}
        \tmr(\dij) & = \tmr(D_1)+\tmr(D_2)
        \\ & = \tmr(D_1^-)+\tmr(D_2^-)
        \\ & = \tmr(D_1^-\rightarrow D_2^-)
        \\ & = \tmr(D_2\rightarrow D_1).
    \end{align*}
\end{proof}

\begin{proof}[Proof of~\Cref{lem:tmr sum}]
    We proceed by induction on $t$. The case $t=2$ holds by \Cref{lem: tmr + tmr}, \Cref{cor:tmr(1->k)} and the proposition. 

    Since either $\tmr(D_1)\leq2$ or $\tmr(D_t)\leq2$, the cases $t>2$ hold by the proposition and the inductive hypothesis.
\end{proof}

\section{Inversion Number for Tournaments}\label{sec:inv tourney}

In this section we prove a version of Theorem~\ref{thm:dijoin} for tournaments, and introduce a couple of lemmas that will be used in the proof of~\Cref{thm:k-join}.

\begin{thm}[\Cref{thm:dijoin} for Tournaments]\label{thm:2 tournaments}
    Let $D_1$ and $D_2$ be tournaments such that $\inv(D_1)=2$ and $\inv(D_2)=k$, then either: 
    \begin{enumerate}[i.]
        \item $\inv(\dij)=k+2$; or
        \item $\inv(\dij)=k+1$.
    \end{enumerate}
    Moreover, $\inv(\dij)=k+1$ if and only if $D_2$ is overfull, and so $k$ is odd. 
\end{thm}

The theorem follows directly from the following result.

\begin{lem}\label{lem: tourney dijoin}
    Let $D_1$ and $D_2$ be tournaments such that $\inv(D_1)=2$ and $\inv(D_2)=k$, then $\inv(\dij)=\tmr(D_2)+2$.
\end{lem}
\begin{proof}
    By \Cref{cor:BJMO main}, we have $\tmr(D_1)=2$. By \Cref{lem: tmr + tmr}, $\tmr(\dij)=\tmr(D_2)+2$. Therefore, $\inv(\dij)\geq\tmr(D_2)+2$ and by \Cref{lem: inv + tmr}, this must be an equality.
\end{proof}

\begin{proof}[Proof of \Cref{thm:2 tournaments}]
By \Cref{lem: tourney dijoin}, we have $\inv(\dij) = \tmr(D_2) +2$. Now, by  \Cref{cor:BJMO main}, $\tmr(D_2)=k$ or $k-1$. Moreover, by \Cref{cor: YWL to tmr},  $\tmr(D_2)=k-1$ if and only if $D_2$ is overfull, in which case $k$ is odd. 
\end{proof}
We can combine \Cref{lem: tourney dijoin} with \Cref{cor: YWL to tmr} and \Cref{lem: wangb}.

\begin{cor}\label{cor:equivalence}
    Let $D_1,D_2$ and $D$ be tournaments such that $\inv(D_1)=1$ and $\inv(D_2)=2$. Then the following are equivalent:
    \begin{enumerate}
        \item $D$ is overfull
        \item $\inv(D)=\tmr(D)+1$
        \item $\inv(D_1\rightarrow D)=\inv(D)$
        \item $\inv(D_2\rightarrow D)=\inv(D)+1$
    \end{enumerate}
\end{cor}

In order to prove \Cref{thm:k-join}, the following two simple Lemmas will be useful. 



\begin{lem}\label{lem: tourney switch}
    Let $D_1$ and $D_2$ be tournaments such that $\inv(D_1)=1\text{ or }2$. Then $\inv(\dij)=\inv(D_2\rightarrow D_1)$.
\end{lem}

\begin{proof}[Proof of \Cref{lem: tourney switch}]
    \Cref{lem: wangb} shows that, if $\inv(D_1)=1$, then $\inv(\dij)=\inv(D_2\rightarrow D_1)$, for any $D_2$. We prove the same is true for $\inv(D_1)=2$.
    
    Denote by $D^-$ the tournament obtained by reversing all arcs of a tournament $D$. Notice that $D$ and $D^-$ have the same set of decycling families and decycling matrices, so $\inv(D)=\inv(D^-)$ and $\tmr(D)=\tmr(D^-)$.
    Then by \Cref{lem: tourney dijoin}:
    \begin{align*}
        \inv(\dij) & = 2+\tmr(D_2)
        \\ & = 2+\tmr(D_2^-)
        \\ & = \inv(D_1^-)+\tmr(D_2^-)
        \\ & = \inv(D_1^-\rightarrow D_2^-)
        \\ & = \inv(D_2\rightarrow D_1)
    \end{align*}
\end{proof}

\begin{lem}\label{lem: tourney 3-join}
    Let $D_1,D_2$ be tournaments with $\inv(D_1),\inv(D_2)\in\{1,2\}$. Then, for any tournament $D_3$, $\inv([D_1,D_2,D_3])=\inv(D_1)+\inv(D_2\rightarrow D_3)$.
\end{lem}





\begin{proof}[Proof of \Cref{lem: tourney 3-join}]
        Clearly $\inv([D_1,D_2,D_3])\leq\inv(D_1)+\inv(D_2\rightarrow D_3)$. For the converse, notice that 
        \begin{align*}
            \inv([D_1,D_2,D_3])& \geq\tmr([D_1,D_2,D_3])
            \\ & = \tmr(D_1) + \tmr(D_2) + \tmr(D_3) && \text{by \Cref{lem:tmr sum}}
            \\ & = \inv(D_1) + \inv(D_2) + \tmr(D_3) && \text{by \Cref{lem: small_inv_not_overfull} \& \Cref{cor: YWL to tmr}}
            \\ & \geq \inv(D_1) + \inv(D_2\rightarrow D_3) && \text{by \Cref{lem: inv + tmr}}
        \end{align*}
\end{proof}

\section{Inversion Number for Oriented Graphs}\label{sec:inv all}

Our main goal in this section is to prove our main theorems for all oriented graphs. 
In order to prove these theorems, we use two facts about the inversion number. 

\begin{obs}\label{obs:subgraph}
    If $D_1\subseteq D_2$, then $\inv(D_1)\leq\inv(D_2)$.
\end{obs}

\begin{prop}\label{prop:extend to tourney}
    For every oriented graph $D$, there is a tournament $D^*$ on the same vertex set, with $D\subseteq D^*$ and $\inv(D^*)=\inv(D)$. 
\end{prop}
\begin{proof}
    Let $\inv(D)=m$, with a decycling family $X_1, X_2, \dots,X_m$. Let $E$ be the acyclic oriented graph obtained by inverting this family. Then there exists a transitive tournament $T$, also on $V(D)$, that contains $E$. If we invert $X_1,X_2,...,X_m$ in $T$, we will obtain a tournament, $D^*$, containing $D$. By \Cref{obs:subgraph}, $\inv(D^*)\geq\inv(D)$, and we already have an $m$-decycling family, so $\inv(D^*)=\inv(D)$. 
\end{proof}

First we prove \Cref{thm:dijoin}, restated below for convenience.

\dijoin* 

\begin{proof}
The case where $\inv(D_1) = 1$ is covered by \Cref{thm:Wang k-join}. Therefore, we will assume that $\inv(D_1) = 2$. 

Set $k = \inv(D_2)$. 
    Clearly $\inv(\dij)\leq k+2$. We first show that $\inv(\dij)\geq k+1$. 
    
    By \Cref{prop:extend to tourney}, we can extend $\dij$ to a tournament $(\dij)^*$, with the same inversion number. Since all the arcs from $V(D_1)$ to $V(D_2)$ remain unchanged, we must have that $(\dij)^*=E_1\rightarrow E_2$, where $E_1$ and $E_2$ are tournaments containing $D_1$ and $D_2$ respectively.

    By \Cref{obs:subgraph}, $\inv(E_1)\geq2$, and so by \Cref{cor:BJMO main}, $\tmr(E_1)\geq2$. Similarly, $\inv(E_2)\geq k$ and so by \Cref{cor:BJMO main}, $\tmr(E_2)\geq k-1$. Therefore, by \Cref{cor:tmr bound},
    \[
    k+1\leq\tmr((\dij)^*)\leq\inv((\dij)^*)=\inv(\dij)
    \]
    as required. 
    
    Now, suppose $\inv(\dij)=k+1$. Then $\inv(E_1\rightarrow E_2)=k+1$. Since $\tmr(E_1) \ge 2$, by \Cref{cor:tmr bound} we know that \[\tmr(E_2) \le \tmr(E_1 \rightarrow E_2) - 2 \le \inv(E_1 \rightarrow E_2) - 2 = k-1.\] Therefore, $\tmr(E_2) = k-1$ which in turn implies $\inv(E_2) = k$. By \Cref{cor:equivalence}, we conclude that $\inv(\overrightarrow{C_3}\rightarrow E_2)=k$ and so $\inv(\overrightarrow{C_3}\rightarrow D_2)=k$ and $k$ is odd.

    Conversely, suppose $\inv(\overrightarrow{C_3}\rightarrow D_2)=k$. By \Cref{prop:extend to tourney}, we can extend $\overrightarrow{C_3}\rightarrow D_2$ to a tournament $(\overrightarrow{C_3}\rightarrow D_2)^*$, this is clearly of the form $\overrightarrow{C_3}\rightarrow E_2$, where $E_2$ is a tournament containing $D_2$. Let $D_1^*$ be a tournament containing $D_1$, with $\inv(D_1^*)=\inv(D_1)=2$. By \Cref{cor:equivalence}, we get $\inv(D_1^*\rightarrow E_2)=k+1$ and so $\inv(\dij)=k+1$.
\end{proof}

We would now like to go on and prove \Cref{thm:k-join}. 

\kjoin*

We first require two lemmas. 

\begin{lem}\label{lem: OG switch}
    Let $D_1$ and $D_2$ be oriented graphs such that $\inv(D_1)=1\text{ or }2$. Then $\inv(\dij)=\inv(D_2\rightarrow D_1)$.
\end{lem}

\begin{lem}\label{lem: OG 3-join}
    Let $D_1,D_2$ be oriented graphs with $\inv(D_1),\inv(D_2)\in\{1,2\}$. Then, for any oriented graph $D_3$, $\inv([D_1,D_2,D_3])=\inv(D_1)+\inv(D_2\rightarrow D_3)$.
\end{lem}

Before proving these lemmas, we show that they give us enough to prove \Cref{thm:k-join}

\begin{proof}[Proof of \Cref{thm:k-join}]
    We can assume without loss of generality that none of the oriented graphs $D_i$, $i \ne j$, have $\inv(D_i)=0$. If some do, then whenever we see adjacent oriented graphs with $\inv(D_i)=0$ and $\inv(D_{i+1})>0$ (or vice-versa), we just treat $D_i\rightarrow D_{i+1}$ as the single oriented graph, with inversion number $\inv(D_{i+1})$. If $\inv(D_i)=0$ for all $i \ne j$, then the theorem trivially holds.

    We proceed by induction on $t$. 
    
    The case $t=2$ holds, by \Cref{thm:dijoin} and \Cref{lem: OG switch}.

    For $t\geq3$, we can assume $j\neq1\text{ or }t$. (If $j=1$, then by \Cref{lem: OG switch}, we consider $[D_t,D_1,D_2,\dots,D_{t-1}]$ and if $j=t$ we consider $[D_2,D_3,\dots,D_t,D_1]$).  

    Let $D=[D_2,D_3,\dots,D_{t-1}]$. Then, \begin{align*}
    \inv([D_1,\dots,D_t])&=\inv([D_1,D,D_t])
    \\ & =\inv([D_t,D_1,D]) && \text{by \Cref{lem: OG switch}}
    \\ & = \inv(D_t)+\inv([D_1,\dots,D_{t-1}]) && \text{by \Cref{lem: OG 3-join}}
    \\ & = \inv(D_t)+\sum_{i=1}^{n-1}\inv(D_i) - \1(D_j \text{ is overfull}) && \text{by induction.}
    \end{align*}

\end{proof}

The case where $\inv(D_j)\leq2$ proves \Cref{conj:alon}, since by~\Cref{lem: small_inv_not_overfull} such $D_j$ are not overfull.

We now prove Lemmas~\ref{lem: OG switch} and \ref{lem: OG 3-join}.

\begin{proof}[Proof of \Cref{lem: OG switch}]
    By \Cref{prop:extend to tourney}, we can extend $\dij$ to a tournament $(\dij)^*$, with the same inversion number. Since the arcs from $D_1$ to $D_2$ remain unchanged, we have $(\dij)^*=E_1\rightarrow E_2$, where $E_1$ and $E_2$ are tournaments containing $D_1$ and $D_2$ respectively. We can also find a tournament $D_1^*$, containing $D_1$, with $\inv(D_1^*)=\inv(D_1)$. We then get 
    \begin{align*}
        \inv(\dij) & = \inv((\dij)^*) \\
        & = \inv(E_1\rightarrow E_2) \\
        & \geq \inv(D_1\rightarrow E_2) \\
        & =\inv(D_1) + \tmr(E_2) && \text{by \Cref{thm:dijoin} \& \Cref{cor: YWL to tmr}}\\
        & = \inv(D_1^*)+\tmr(E_2) \\
        & = \inv(D_1^*\rightarrow E_2) \\
        & = \inv(E_2\rightarrow D_1^*) && \text{by \Cref{lem: tourney switch}}\\
        & \geq \inv(D_2\rightarrow D_1). 
    \end{align*}
    Notice that $\inv(D_2\rightarrow D_1)=\inv(D_1^-\rightarrow D_2^-)$, and so applying the same argument we get the reverse inequality:
    \[
    \inv(D_2\rightarrow D_1)=\inv(D_1^-\rightarrow D_2^-)\geq\inv(D_2^-\rightarrow D_1^-)=\inv(\dij).
    \]
\end{proof}

\begin{proof}[Proof of \Cref{lem: OG 3-join}] 
    Clearly $\inv([D_1,D_2,D_3])\leq\inv(D_1)+\inv(D_2\rightarrow D_3)$.

    For the converse, by \Cref{prop:extend to tourney}, we can extend $[D_1,D_2,D_3]$ to a tournament $[D_1,D_2,D_3]^*$, with the same inversion number. Clearly $[D_1,D_2,D_3]^*=[E_1,E_2,E_3]$, where $E_1, E_2$ and $E_3$ are tournaments containing $D_1, D_2$ and $D_3$ respectively.

    We can also find a tournament $D_1^*$, containing $D_1$, with $\inv(D_1^*)=\inv(D_1)$. Then we get:
    \begin{align*}
        \inv([D_1,D_2,D_3]) & = \inv([D_1,D_2,D_3]^*)
        \\ & = \inv([E_1,E_2,E_3])
        \\ & \geq \inv([D_1,E_2,E_3])
        \\ & = \inv(D_1) + \tmr(E_2\rightarrow E_3) && \text{by \Cref{thm:dijoin} \& \Cref{cor: YWL to tmr}}
        \\ & = \inv(D_1^*) + \tmr(E_2\rightarrow E_3)
        \\ & = \inv([D_1^*,E_2,E_3])
        \\ & = \inv(D_1^*) + \inv(E_2\rightarrow E_3) && \text{by \Cref{lem: tourney 3-join}}
        \\ & \geq \inv(D_1) + \inv(D_2\rightarrow D_3).
    \end{align*}
\end{proof}

\section{Further Ideas and Open Questions}\label{sec:Conjectures}

In this section we discuss some conjectures and give initial thoughts on proving them. 

It is now known that the dijoin conjecture, $\inv(\dij)=\inv(D_1)+\inv(D_2)$, holds for the following pairs ($\inv(D_1),\inv(D_2)$):
\[
(0,k),(k,0),(1,1),(1,2k),(2k,1),(2,2k),(2k,2)
\]

We also know, from \citep{alon2024} and \citep{aubian2025}, that there are counterexamples whenever $\inv(D_1)$ or $\inv(D_2)$ is odd $\geq3$. This leaves the case $(2j,2k)$. One approach to prove the conjecture in this case would be to show that the additivity property of tournament minimum rank holds generally.

\begin{conj}\label{conj: ours}
    For all tournaments $D_1$, $D_2$, we have $\tmr(\dij) = \tmr(D_1)+\tmr(D_2).$
\end{conj}

The $(2j,2k)$ case of the dijoin conjecture would then follow from \Cref{cor:BJMO main}.

To prove \Cref{conj: ours}, we are able to follow our proof of the additivity property for $\inv(D_1)=2$ up until the final step: verifying that $B'$ is still a decycling matrix for $D_2$.

Firstly, note that a more general version of \Cref{lem: submatrix} holds. Let $A$ be an $n\times n$ symmetric matrix with entries in $\F_2$, and rank $r$. Then it is always possible  to find an $r\times r$ principal submatrix $A'$ of rank $r$ -- see e.g.~\citep{barrett2014}.

    


We can then sketch our potential proof for \Cref{conj: ours}.

\begin{itemize}
    \item Suppose $\tmr(\dij)<\tmr(D_1)+\tmr(D_2)$, for a contradiction.
    \item Then take a decycling matrix, $M=\begin{psmallmatrix}
        A & C \\ {C^T} & B
    \end{psmallmatrix}$, for $\dij$ of rank ${<\tmr(D_1)+\tmr(D_2)}$. 
    \item Let $A'$ be a full rank principal submatrix of $A$, with rank $\geq\tmr(D_1)$, and apply row/column additions to get a matrix $\begin{psmallmatrix}
        A' & 0 \\ 0 & B'
    \end{psmallmatrix}$, still with rank $<\tmr(D_1)+\tmr(D_2)$.  
    \item Here $B'$ will have rank at most $\tmr(D_2)-1$, so if it is still a decycling matrix for $D_2$, then we would have a contradiction.
\end{itemize}

Unfortunately, it is not straightforward to show this final step, so a further investigation of $B'$ is required.

Notice that the column/row operations on $M$ amount to a transformation $P^TMP$, where $P$ is of the form $\begin{psmallmatrix}
    I & X \\ 0 & I
\end{psmallmatrix}$. In particular we get:
\[
\begin{pmatrix}
    I&0\\X^T&I
\end{pmatrix}
\begin{pmatrix}
    A'&C\\C^T&B
\end{pmatrix}
\begin{pmatrix}
    I&X\\0&I
\end{pmatrix}
=\begin{pmatrix}
    A'&A'X+C\\X^TA'+C^T&X^TA'X+X^TC+C^TX+B
\end{pmatrix}
\]
Since we eliminate the off-diagonal blocks, we are using $X=A'^{-1}C$, so our final matrix is 
\[
\begin{pmatrix}
    A' & 0 \\ 0 & C^TA'^{-1}C+B
\end{pmatrix}.
\]
So we are interested in showing that adding $C^TA'^{-1}C$ takes the decycling matrix $B$ to another decycling matrix. This quantity is quite mysterious, in particular some insight into the combinatorial meaning of the inverse or product of decycling matrices would be useful. 

It's worth noting that a case by case check of $3\times3$ submatrices $A'$ gives us that $C^TA'^{-1}C + B$ is no longer a decycling matrix for $D_2$, if and only if $A'$ is a decycling matrix for $\overrightarrow{C_3}$. 

It is not obvious whether this property holds in general, though. Specifically, is $C^TA'^{-1}C + B$ not a decycling matrix for $D_2$ if and only if $A'$ is a decycling matrix for some non-transitive tournament?

A proof of \Cref{conj: ours} would also give us the lower bound, $\inv(\dij)\geq \inv(D_1)+\inv(D_2)-2$. Notice that this lower bound would only possibly be attained if $\inv(D_1)=\tmr(D_1)+1$ and $\inv(D_2)=\tmr(D_2)+1$. We conjecture that it is never attained: 

\begin{conj}
    Let $D_1$ and $D_2$ be tournaments, then
    \[
    \inv(\dij)\geq\inv(D_1)+\inv(D_2)-1
    \]
    with equality if and only if $\inv(D_1)=\tmr(D_1)+1$ or $\inv(D_2)=\tmr(D_2)+1$.
\end{conj}

The reason we suspect this to be the case is that, if $\inv(\dij)=\inv(D_1)+\inv(D_2)-2$, then it could not have a decycling matrix of the form $\begin{psmallmatrix}
    A&0\\0&B
\end{psmallmatrix}$, where $A$ and $B$ are minimum rank decycling matrices for $D_1$ and $D_2$ respectively. This is because $A$ and $B$ would both have all zero diagonals contradicting $\inv(\dij)=\tmr(\dij)$.

\section{Acknowledgements}
 While conducting this research, the first author was supported by the European Research Council (ERC) under the European Union Horizon 2020 research and innovation programme (grant agreement No.\ 947978).

 The second author was supported by the Warwick Mathematics Institute through the Undergraduate Research Support Scheme (URSS).

\bibliography{main}
\end{document}